\documentclass[12pt]{article}
\usepackage[a4paper]{geometry}
\usepackage[italian,english]{babel}

\usepackage{amsmath, accents}
\usepackage{amsfonts}
\usepackage{amstext}
\usepackage{amssymb}
\usepackage{amsthm}
\usepackage{amscd}
\usepackage{mathrsfs}

\usepackage[pagebackref,draft=false]{hyperref}
\hypersetup{colorlinks,
linkcolor=myrefcolor,
citecolor=mycitecolor,
urlcolor=myurlcolor}

\usepackage[capitalize]{cleveref}
\usepackage{caption}
\usepackage{etaremune}

\usepackage{xcolor}
\definecolor{myurlcolor}{rgb}{0,0,0.4}
\definecolor{mycitecolor}{rgb}{0,0.5,0}
\definecolor{myrefcolor}{rgb}{0.5,0,0}
\usepackage{graphicx}
\usepackage{tikz}
\usepackage{tikz-cd}

\usepackage{etoolbox}
\usepackage{makeidx}
\usepackage{sectsty}
\usepackage{dsfont}
\usepackage{enumitem} 
\usepackage[]{latexsym}
\usepackage{braket}
\usepackage{caption}
\usepackage[utf8]{inputenx}
\usepackage[T1]{fontenc}
\usepackage{lmodern}
\usepackage{textcomp}
\usepackage{microtype}
\usepackage{totcount}
\usepackage{blindtext}

\newtheorem{remark}{Remark}

\newtheorem{theorem}{Theorem}

\newtheorem{proposition}{Proposition}

\newtheorem*{proof*}{Proof}

\newcommand{\be}{\begin{equation}}
\newcommand{\ee}{\end{equation}}
\newcommand{\bea}{\begin{eqnarray}}
\newcommand{\eea}{\end{eqnarray}}

\newcommand{\m}{\mathfrak{m}}

\newcommand{\lag}{\mathscr{L}}

\newcommand{\dd}{{\rm d}}
\newcommand{\de}{\partial}

\title{The inverse problem for a class of implicit differential equations and the coisotropic embedding theorem}

\author{L. Schiavone$^{1,2}$ \href{https://orcid.org/0000-0002-1817-5752}{\includegraphics[scale=0.7]{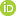}}, \\ 
\footnotesize{$^{3}$\textit{ Departamento de Matem\'aticas, Universidad Carlos III de Madrid, Legan\'es, Madrid, Spain}}  \\
\footnotesize{$^{2}$\textit{ e-mail: \texttt{lschiavo[at]math.uc3m.es}}} 
}

\begin{document}

\maketitle

\begin{abstract}
We carry on the approach used in \cite{Schiavone2023-Inverse_problem_Electrodynamics} to provide a solution for the inverse problem of the calculus of variations for Maxwell equations in vacuum and we provide an abstract theory including all implicit differential equations that can be formulated in terms of vector fields over pre-symplectic manifolds.
\end{abstract}

\tableofcontents

\section*{Introduction}
\label{Sec:Introduction}
\addcontentsline{toc}{section}{\nameref{Sec:Introduction}}

In a previous paper \cite{Schiavone2023-Inverse_problem_Electrodynamics} the author presented the coisotropic embedding theorem as a tool to provide a solution for the inverse problem of the calculus of variations for the system of Maxwell equations in vacuum.
As recalled in the introduction of \cite{Schiavone2023-Inverse_problem_Electrodynamics}, the inverse problem of the calculus of variations is a very old problem in mathematical analysis and there is a substantial body of literature on its solution.
We refer to the introduction of \cite{Schiavone2023-Inverse_problem_Electrodynamics} for a relatively comprehensive list of references and for some historical notes about the subject.
Here we limit ourselves to recall that a solution for the inverse problem of the calculus of variations for a quite general class of partial differential equations exists and is elegantly resumed in \textit{D. Krupka}'s book on \textit{Global Variational Geometry} \cite{Krupka2015-Variational_Geometry} (see also the previous remarkable results about the inverse problem in \cite{Helmholtz1887-Inverse_Problem, Douglas1941-Inverse_Problem, Henneaux1982-Inverse_Problem, Crampin1981-Inverse_Problem, Morandi-Ferrario-LoVecchio-Marmo-Rubano1990-Inverse_Problem, Ibort-Solano1991-Inverse_Problem_coupled_systems}).
It takes advantage of the fact that the space of solutions of a large class of PDEs of order $k$ can be formulated as the set of zeros of a family of smooth functions, say $E_\sigma$ of the type:
\be \label{Eq:k-order pde}
E_\lambda(x^j,\, y^\sigma,\, y^\sigma_{j_1},\, y^\sigma_{j_1 j_2},\, ...,\, y^\sigma_{j_1...j_k}) \,,
\ee
where the $x^j$ ($j=1,..., n$) are the independent variables of the differential equations, whereas $y^\sigma$ ($\sigma=1,...,l$) are the dependent ones and $y^\sigma_{j_1...j_k}$ are the $k$-th order derivatives of the $y^\sigma$ with respect to the $x^j$.
The set of the independent variables, the dependent ones and their derivatives up to order $k$ can be seen as a set of fibered coordinates over the jet bundle of order $k$, say $\mathbf{J}^k \pi$, of the fibre bundle $\pi \, :\; Y \to X$ with $\left\{\,x^j,\, y^\sigma \right\}_{j=1,...,n;\,\sigma=1,...,l}$ being fibered coordinates on $Y$ associated to the set of local coordinates $\left\{\, x^j \,\right\}_{j=1,...,n}$ on $X$. 
The map $\pi$ above is the projection onto the first factor.
The functions $E_\lambda$ can be seen as coefficients of the following $(n+1)$-form on $\mathbf{J}^k \pi$:
\be
E \,=\, E_\lambda \eta^\lambda \wedge \underbrace{\dd x^1 \wedge ... \wedge \dd x^n}_{vol_X} \,,
\ee
usually referred to as a \textit{source form} on $\mathbf{J}^k \pi$, where $\eta^\lambda \,=\, \dd y^\lambda - y^\lambda_j \dd x^j$.
In this setting it is proven that under suitable assumptions on the differential form $E$, a Lagrangian whose extrema are the solutions of the PDE considered above exists \cite[Theorem 13]{Krupka2015-Variational_Geometry}.

\noindent In \cite{Schiavone2023-Inverse_problem_Electrodynamics}, we took a completely different approach and we took advantage of the fact that solutions of Maxwell equations in vacuum can be seen as the space of integral curves of a family of vector fields on a pre-symplectic manifold.
In particular we recalled how a solution of the inverse problem of the calculus of variations for first order ordinary differential equations exists and it does not require to model the PDE as the coefficients of a source form on a jet bundle.
It is resumed in the following theorem \cite{Ibort-Solano1991-Inverse_Problem_coupled_systems, Ciglia-DC-Ibort-Marmo-Schiav-Zamp2020-Heisenberg}.
Recall that a first order ordinary differential equation over a smooth manifold is represented by a vector field whose integral curves are the solutions of the equation.
\begin{theorem}[\textsc{Lagrangian for first order odes}]\label{Thm:Lagrangian for first order odes}
Given a vector field $\Gamma$ over a manifold $\mathcal{M}$, if there exists a symplectic form $\omega$ on $\mathcal{M}$ such that:
        \be
        \mathcal{L}_\Gamma \omega \,=\, 0 \,,
        \ee
        then, there exists a local Lagrangian function $\lag$ on $\mathbf{T}\mathcal{M}$ such that the pre-symplectic Hamiltonian system $(\mathbf{T}\mathcal{M},\, \omega_\lag,\, E_\lag)$ yields, by applying the pre-symplectic constraint algorithm (see \cite{Gotay-Nester-Hinds1978-DiracBergmann_constraints} for the definition of the algorithm), that:
        \begin{itemize}
            \item $\mathcal{M}_\infty \equiv \mathcal{M}$ is the final stable manifold of the algorithm;
            \item ${\omega_\lag}_\infty \,=\, \mathfrak{i}_\infty^\star \omega_\lag \equiv \omega$, where $\mathfrak{i}_\infty$ is the immersion of $\mathcal{M}_\infty \,\equiv \, \mathcal{M}$ into the original manifold $\mathbf{T}\mathcal{M}$;
            \item $\Gamma$ satisfies $i_\Gamma {\omega_\lag}_\infty \,=\, i_\Gamma \omega \,=\, \dd {E_\lag}_\infty$ where ${E_\lag}_\infty \,=\, \mathfrak{i}_\infty^\star E_\lag$.
        \end{itemize}
        This means that, given a vector field $\Gamma$ over a manifold, if the manifold is symplectic and if the vector field preserves the symplectic structure, there exists a degenerate Lagrangian on the tangent bundle of the symplectic manifold giving rise to first order Euler-Lagrange equations coinciding with the set of ordinary differential equations defining $\Gamma$.
\end{theorem}

\noindent With the latter theorem in hands, we showed how, starting from the family of vector fields associated to the system of Maxwell equations, one may use the coisotropic embedding theorem to obtain a vector field enconding the same set of differential equations over a suitable (larger) symplectic manifold fulfilling the hypothesis of \cref{Thm:Lagrangian for first order odes} and for which a Lagrangian description can be given.

The aim of the present paper it to push forward this point of view to a more abstract setting that does not only include the single example considered in \cite{Schiavone2023-Inverse_problem_Electrodynamics}. 
In particular we will prove that, via a procedure with the same philosophy, a solution for the inverse problem of the calculus of variations can be given for a class of implicit partial differential equations which is interesting from the physical point of view since it includes all partial differential equations emerging as equations of motion within classical field theories associated with the quantum field theories describing fundamental interactions.

The structure of the paper is the following.

\noindent In \cref{Sec:The coisotropic embedding theorem} we recall the content of the coisotropic embedding theorem.
In particular, in \cref{Subsec:The symplectic thickening} we recall how to construct the so called \textit{symplectic thickening} of a pre-symplectic manifold and in \cref{Subsec:Lift of vector fields to the symplectic thickening} we show how a class of vector fields over the pre-symplectic manifold can be canonically lifted to the symplectic thickening.

\noindent In \cref{Sec:The inverse problem for a class of implicit differential equations} we prove the main result of the paper, namely how to solve the inverse problem of the calculus of variations for all those implicit partial differential equations associated to pre-symplectic Hamiltonian systems.

\section{The coisotropic embedding theorem}
\label{Sec:The coisotropic embedding theorem}

In this section, we recall the main tool used along the manuscript to provide a solution for the inverse problem, i.\ e. the \textit{coisotropic embedding theorem}.

\noindent In particular, in \cref{Subsec:The symplectic thickening} we recall how to embed a pre-symplectic manifold into a larger symplectic one referred to as its \textit{symplectic thickening} and in \cref{Subsec:Lift of vector fields to the symplectic thickening} we recall how to lift a particular class of vector fields from the pre-symplectic manifold to the symplectic thickening.

\subsection{The symplectic thickening}
\label{Subsec:The symplectic thickening}

The original version of the \textit{coisotropic embedding theorem} is due to \textit{M. Gotay} and can be found in \cite{Gotay1982-Coisotropic_embedding}.
Other versions have appeared in the literature in the subsequent years.
An interesting equivariant version of the theorem can be found in \cite{Guillemin-Sternberg1990-Symplectic_techniques}.
Another, slightly more modern, version can be found in \cite{Oh-Park2005-Embedding_coisotropo} and its application to the construction of Poisson brackets on the space of solutions of the equations of motion within gauge theories can be found in the contributions of the author \cite{Ciaglia-DC-Ibort-Marmo-Schiav-Zamp2021-Cov_brackets_toappear, Ciaglia-DC-Ibort-Mar-Schiav-Zamp2022-Non_abelian, Ciaglia-DC-Ibort-Mar-Schiav-Zamp2022-Palatini}.
The author also used the same technique in \cite{Ciaglia-DC-Ibort-Marmo-Schiav-Zamp-2022-Symmetry} to obtain a Noether type theorem for pre-symplectic Hamiltonian systems.

\begin{theorem}[\textsc{Coisotropic embedding theorem}]
Let us consider a smooth pre-symplectic (Banach) manifold $(\mathcal{M},\, \omega)$.
Let $K_m \,=\, \mathrm{ker} \, \omega$ be isomorphic at each $m \in \mathcal{M}$.
Consider, at every $m \in \mathcal{M}$, a complement of $K_m$ into $\mathbf{T}_m \mathcal{M}$, say $W_m$, and let it be the kernel of a projection map:
\be
P \;\; :\;\; \mathbf{T}\mathcal{M} \to K \,.
\ee
Let $\pi \;:\;\; \mathbf{K} \to \mathcal{M}$ denote the characteristic bundle over $\mathcal{M}$, i. e.\ the vector bundle over $\mathcal{M}$ having $K_m$ as typical fibre.
Denote by $\tau \;:\;\; \mathbf{K}^\star \to \mathcal{M}$ the dual bundle to $\pi$.

\noindent Then, there exists a neighborhood of the zero section of $\tau$, that we denote by $\tilde{\mathcal{M}}$, upon which the symplectic form below can be defined:
\be
\tilde{\omega} \,=\, \tau^\star \omega + \dd \vartheta^P \,,
\ee
with:
\be
\vartheta^P_\mu (X) \,=\, \mu\left(\, P \circ T\tau (X) \,\right) \,,
\ee
for $\mu \in \mathbf{K}^\star$, $X \in \mathbf{T}_\beta \mathbf{K}^\star$ and $T\tau$ being the tangent map of $\tau$:
\be
T\tau \;\; :\;\; \mathbf{T}\mathbf{K}^\star \to \mathbf{T}\mathcal{M} \,.
\ee
The smooth (Banach) manifold $(\tilde{\mathcal{M}},\, \tilde{\omega})$ is a (generally weakly) symplectic manifold called the \textbf{symplectic thickening} of $(\mathcal{M},\, \omega)$) for which $(\mathcal{M},\, \omega)$ is a coisotropic submanifold, namely $\tau^\star \tilde{\omega} \,=\, \omega$.
\end{theorem}
In a system of local coordinates the projection $P$ can be expressed in terms of an idempotent $(1,\,1)$-tensor field on $\mathcal{M}$ of the type:
\be
P \,=\, P_j \otimes V^j \,,
\ee
where $\left\{\, V^j \,\right\}_{j=\mathcal{J}}$\footnote{$\mathcal{J}$ denotes a suitable (eventually infinite-dimensional) index set.} is a basis for $K_m$ for all $m \in \mathcal{M}$ and $\left\{\, P_j \,\right\}_{j=\mathcal{J}}$ is a collection of $1$-forms on $\mathcal{M}$ providing the components along $K$ associated to the chosen connection in the above mentioned basis of any vector field $Y$:
\be
P(Y) \,=\, P_j(Y) V^j \,=:\, {Y^v}_j V^j \,. 
\ee
As it was shown in \cite{Ciaglia-DC-Ibort-Mar-Schiav-Zamp2022-Non_abelian, Schiavone2023-Inverse_problem_Electrodynamics} $\tilde{\omega}$ locally reads:
\be
\tilde{\omega} \,=\, \tau^\star \omega + \dd \mu^j \wedge P_j + \mu^j \dd P_j \,,
\ee
where $\left\{\, \mu_j \,\right\}_{j \in \mathscr{J}}$ is a suitable basis for $K_m^\star$.

\subsection{Lift of vector fields to the symplectic thickening}
\label{Subsec:Lift of vector fields to the symplectic thickening}

In this section we define a way of lifting vector fields from a pre-symplectic manifold $(\mathcal{M},\, \omega)$ to its symplectic thickening $(\tilde{\mathcal{M}},\, \tilde{\omega})$ that will be useful along the manuscript.

\noindent Consider a vector field $X \in \mathfrak{X}(\mathcal{M})$ and denote by $P$ the connection used in constructing the symplectic thickening. 
By virtue of $P$, $X$ splits into a vertical and a horizontal part, say:
\be
X \,=\, X_V + X_H \,,
\ee
where:
\be
P(X_V) \,=\, X_V \,, \;\; \qquad \;\; P(X_H) \,=\, 0 \,. 
\ee
The differential $1$-form $\vartheta^P$ can be used to lift the vertical part $X_V$ to $\tilde{\mathcal{M}}$ to the vector field $\tilde{X}_V$ defined by:
\be \label{Eq:lift vertical part}
\begin{split}
T \tau (\tilde{X}_V) \,&=\, X_V \,, \\
\mathcal{L}_{\tilde{X}_V} \vartheta^P \,&=\, 0 \,.
\end{split}
\ee
This lift is well defined since, using Cartan's formula for the Lie derivative and the definition of $\vartheta^P$, the latter equation reads:
\be \label{Eq:condition lift vertical part}
\tilde{X}_V \left[\, \mu(P(Y_V)) \,\right] \,=\, \mu \left[\, 
P([X_V,\,Y]) \,\right] \;\;\; \forall \,\, Y \in \mathfrak{X}(\mathcal{M}) \,.
\ee
The right hand side of the latter equation is known and, since the equality must hold for any vector field $Y$ on $\mathcal{M}$, we can write it for a family of vector fields $\left\{\, Y_j \,\right\}_{j \in \mathscr{J}}$ such that the family $\left\{\, \dd \left[\, 
\mu(P(Y_j)) \,\right] \,\right\}_{j \in \mathscr{J}}$ is a basis of $1$-forms on $K^\star_m$ for each $m$ (where $m \,=\, \tau(\mu)$).
Thus, \eqref{Eq:condition lift vertical part} defines the components of $\tilde{X}_V$ along a basis of differential $1$-forms on the fibres of $\tau$ and, thus, completely defines the lift.

\noindent On the other hand the horizontal part $X_H$ can be lifted to $\tilde{\mathcal{M}}$ to a vector field $\tilde{X}_H$ defined by:
\be \label{Eq:lift horizontal part}
\begin{split}
T \tau (\tilde{X}_H) \,&=\, X_H \,, \\
i_{\tilde{X}_H} \dd \mu_j \,&=\, 0 \;\;\;\;\;\; \forall \,\, j \in \mathscr{J} \,,
\end{split}
\ee
where $\mathscr{J}$ is a suitable index set and $\left\{\, \dd \mu_j \,\right\}_{j \in \mathscr{J}}$ is again a basis of differential $1$-forms on $K^\star_m$ at each $m \in \mathcal{M}$.

\section{The inverse problem for a class of implicit differential equations}
\label{Sec:The inverse problem for a class of implicit differential equations}

Consider dynamical systems whose trajectories are the integral curves of a vector field, say $\Gamma$, over a pre-symplectic manifold $(\mathcal{M},\, \omega)$, satisfying:
\be \label{Eq:Gamma Hamiltonian}
i_{\Gamma} \omega \,=\, \dd \mathscr{H} \,,
\ee
for some $\mathscr{H} \in \mathcal{C}^\infty (\mathcal{M})$.

\noindent As recalled in the concluding section of \cite{Schiavone2023-Inverse_problem_Electrodynamics}, all gauge theories describing fundamental interactions lie in this class of systems.

\noindent To provide a Lagrangian formulation for $\Gamma$ we proceed as follows:
\begin{itemize}
    \item We embed the pre-symplectic manifold $(\mathcal{M},\, \omega)$ into its symplectic thickening $(\tilde{\mathcal{M}},\, \tilde{\omega})$ via the coisotropic embedding theorem described in \cref{Subsec:The symplectic thickening};
    \item We lift the vector field $\Gamma$ to a vector field $\tilde{\Gamma}$ on $\tilde{\mathcal{M}}$ in the way described in \cref{Subsec:Lift of vector fields to the symplectic thickening};
    \item If $\tilde{\Gamma}$ satisfies the condition:
    \be \label{Eq:preservation omegatilde}
    \mathcal{L}_{\tilde{\Gamma}} \tilde{\omega} \,=\, 0 \,,
    \ee
    we will provide a Lagrangian formulation for $\tilde{\Gamma}$ via a direct application of \cref{Thm:Lagrangian for first order odes}.
    We will discuss under which conditions \eqref{Eq:preservation omegatilde} is satisfied.
\end{itemize}

\noindent Regarding the first step, we follow the content of \cref{Sec:The coisotropic embedding theorem} and we construct the symplectic thickening of $(\mathcal{M},\, \omega)$ as the symplectic manifold $(\tilde{\mathcal{M}},\, \tilde{\omega})$, where $\tilde{\mathcal{M}}$ is a suitable tubular neighborhood of the zero section of the dual of the characteristic bundle over $\mathcal{M}$ and $\tilde{\omega}$ reads the symplectic form:
\be
\tilde{\omega} \,=\, \pi^\star \omega + \dd \vartheta^P \,,
\ee
whose local expression is:
\be
\tilde{\omega} \,=\, \pi^\star \omega + \mu^j \dd P_j + \dd \mu^j \wedge P_j \,,
\ee
where $P \,=\, P_j \otimes V^j$ is the connection used along the coisotropic embedding procedure and represented by means of a $1-1$ tensor field on $\mathcal{M}$, $P_j$ are differential $1$-forms, $\left\{\, V^j \,\right\}_{j \in \mathcal{J}}$ is a basis for $\mathrm{ker} \, \omega$ at each point, $\mathcal{J}$ denoting a suitable (eventually infinite-dimensional) index set.

\noindent The second step is a straightforward application of the construction outlined in \cref{Subsec:Lift of vector fields to the symplectic thickening}.
It is important to note that since $\tilde{\Gamma}$ is a lift of $\Gamma$, its integral curves will projects onto the integral curves of $\Gamma$ and, thus, the symplectic manifold $\tilde{\mathcal{M}}$ plays the role of an unfolding space to provide a Lagrangian description for $\Gamma$.

\noindent Regarding the third step, let us first prove the following proposition.
\begin{proposition} \label{Prop:conservation of omega}
The vector field $\tilde{\Gamma}$ constructed above is such that:
\be
\mathcal{L}_{\tilde{\Gamma}} \tilde{\omega} \,=\, 0 \,,
\ee
if $P$ is a flat connection.
\begin{proof}
Taking into account the expression of $\tilde{\omega}$, we have:
\be \label{Eq:Lie derivative omegatilde}
\mathcal{L}_{\tilde{\Gamma}} \tilde{\omega} \,=\, \mathcal{L}_{\tilde{\Gamma}} \left(\,\tau^\star \omega\,\right) + \mathcal{L}_{\tilde{\Gamma}} \dd \vartheta^P \,.
\ee
The first term on the right hand side vanishes since:
\be 
\mathcal{L}_{\tilde{\Gamma}} \left(\,\tau^\star \omega\,\right) \,=\, \tau^\star \left(\, 
\mathcal{L}_{\Gamma} \omega \,\right) \,=\, 0 \,,
\ee
where the first equality is a consequence of:
\be
T\tau \left(\,\tilde{\Gamma}\,\right) \,=\, \Gamma \,,
\ee
and the second equality is a consequence of \eqref{Eq:Gamma Hamiltonian}.

\noindent The second term on the right hand side of \eqref{Eq:Lie derivative omegatilde} splits as:
\be \label{Eq:Lie derivative dtheta}
\mathcal{L}_{\tilde{\Gamma}} \dd \vartheta^P \,=\, \mathcal{L}_{{\tilde{\Gamma}}_V} \dd \vartheta^P + \mathcal{L}_{{\tilde{\Gamma}}_H} \dd \vartheta^P \,.
\ee
The first term on the right hand side of \eqref{Eq:Lie derivative dtheta} reads:
\be
\mathcal{L}_{{\tilde{\Gamma}}_V} \dd \vartheta^P \,=\,  \dd \mathcal{L}_{{\tilde{\Gamma}}_V} \vartheta^P \,=\,0 \,,
\ee
where the last equality is a consequence of \eqref{Eq:lift vertical part}.
The second term on the right hand side of \eqref{Eq:Lie derivative dtheta} reads:
\be \label{Eq:Lie derivative dtheta horizontal}
\mathcal{L}_{{\tilde{\Gamma}}_H} \dd \vartheta^P \,=\, \dd \mathcal{L}_{{\tilde{\Gamma}}_H} \vartheta^P \,=\, \dd i_{{\tilde{\Gamma}}_H} \dd \vartheta^P \,,
\ee
where the differential $1$-form $i_{\tilde{\Gamma}_H} \dd \vartheta^P$ is:
\be \label{Eq: second term lie derivative}
i_{\tilde{\Gamma}_H} \dd \vartheta^P (\tilde{Y}) \,=\, \tilde{\Gamma}_H [\vartheta^P(\tilde{Y})]  - \tilde{Y} [\vartheta^P(\tilde{\Gamma}_H)] - \mu (P([\Gamma_H, Y])) \,.
\ee
The second term on the right hand side of \eqref{Eq: second term lie derivative} vanishes since $\vartheta^P(\tilde{\Gamma}_H) \,=\, \mu(P(\Gamma_H)) \,=\, 0$ since $P(\Gamma_H) \,=\, 0$.
The first term on the right hand side is:
\be \label{Eq: term 1 proof}
\tilde{\Gamma}_H \left[ \vartheta^P(\tilde{Y})\right] \,=\, \tilde{\Gamma}_H \left[ 
\mu(P(Y)) \right] \,=\, \mu \left[ 
\mathcal{L}_{\Gamma_H} (P(Y)) \right] \,,
\ee
where the first equality is a consequence of the definition of $\vartheta^P$ and the second equality is due to \eqref{Eq:lift horizontal part}.
The third term on the right hand side of \eqref{Eq: second term lie derivative} reads:
\be \label{Eq: last term proposition}
-\mu(P([\Gamma_H,\, Y])) \,=\, -\mu(P([\Gamma_H,\, Y_H])) -\mu(P([\Gamma_H,\, Y_V])) \,. 
\ee
The first term on the right hand side of \eqref{Eq: last term proposition} vanishes because, since $P$ is flat, the Lie bracket of two horizontal vector fields is horizontal as well and, thus, $P([\Gamma_H,\, Y_H]) \,=\, 0$.
The second term on the right hand side of \eqref{Eq: last term proposition} reads:
\be \label{Eq: term 2 proof}
-\mu\left(P([\Gamma_H,\,Y_V])\right) \,=\, - \mu \left(\,\mathcal{L}_{\Gamma_H}\left(P(Y)\right)\right) + \mu\left( \mathcal{L}_{\Gamma_H}P(Y)\right) \,=\, - \mu \left(\,\mathcal{L}_{\Gamma_H}\left(P(Y)\right)\right) \, 
\ee
where the first equality is due to the property:
\be
i_{[X,Y]} \theta \,=\, [\mathcal{L}_X,\, i_Y] \theta \,,
\ee
of the Lie derivative and the contraction whereas the second equality is a consequence of \eqref{Eq:lift horizontal part}.

\noindent The thesis follows by noting that the sum of \eqref{Eq: term 1 proof} and \eqref{Eq: term 2 proof} is $0$.
\end{proof}
\end{proposition}

\noindent Some comments have to be made.

\begin{remark}
The condition that $P$ is flat is not the minimal one in order for the second term on the right hand side of \eqref{Eq:Lie derivative dtheta horizontal} to vanish.
Indeed, flatness of $P$ implies that the Lie bracket of any two horizontal vector fields vanishes whereas we only need the Lie bracket of $\Gamma_H$ with any horizontal vector field to vanish.
\end{remark}

\begin{remark} \label{Rem:flat connection}
The existence of a flat connection on a generic fibre bundle is not always guaranteed (see, for instance \cite{Narasimhan1979-SU2, Singer1978-Gribov_ambiguity} for a counterexample). 
However, any trivial bundle admit a flat connection \cite[Example $11.4.1$]{Giac-Mang-Sard2010-Geometric_Classical_Quantum_Mechanics} and any fibre bundle is locally trivial.
Consequently, the bundle $\mathcal{M} \to \mathcal{M}/K$ is locally trivial around any point in the base manifold and if one restricts his attention to such a trivial bundle, the existence of a flat connection is guaranteed.
In that case condition \eqref{Eq:preservation omegatilde} holds only locally on the points of $\tilde{\mathcal{M}}$ whose projections lie on the submanifold of $\mathcal{M}$ salected by restricting to the trivial bundle above.
\end{remark}

\noindent By virtue of \cref{Prop:conservation of omega} and \cref{Rem:flat connection}, the vector field $\tilde{\Gamma}$ on $\tilde{\mathcal{M}}$ satisfies the hypothesis of \cref{Thm:Lagrangian for first order odes}, at least locally on $\tilde{\mathcal{M}}$.
Consequently, there exists a local Lagrangian function $\lag$ on $\mathbf{T}\tilde{\mathcal{M}}$ such that the pre-symplectic Hamiltonian system $(\mathbf{T}\tilde{\mathcal{M}},\, \omega_\lag, E_\lag)$, where:
\be
\omega_\lag \,=\, \dd \dd_S \lag \,, \qquad E_\lag \,=\, \Delta(\lag) - \lag \,,
\ee
with $S$ and $\Delta$ being the soldering one-form and the dilation vector field on $\mathbf{T}\tilde{\mathcal{M}}$, has $\tilde{\mathcal{M}}$ as final stable manifold for the pre-symplectic constraint algorithm, $\mathfrak{j}^\star \omega_\lag \,=\, \tilde{\omega}$ ($\mathfrak{j}$ denoting the immersion of $\tilde{\mathcal{M}}$ into $\mathbf{T}\tilde{\mathcal{M}}$) and $\tilde{\Gamma}$ turns out to be the Hamiltonian vector field associated to $\mathfrak{j}^\star E_\lag$ with respect to the symplectic structure $\tilde{\omega}$.
This means that solutions of Euler-Lagrange equations for $\lag$ are the integral curves of $\tilde{\Gamma}$.

\noindent Finally, since $\tilde{\Gamma}$ is a lift of $\Gamma$, its integral curves projects back onto the integral curves of $\Gamma$.
Therefore, the symplectic manifold $\tilde{\mathcal{M}}$ can be seen as an unfolding space to provide a solution for the inverse problem of the calculus of variations for implicit differential equations that can be formulated in terms of vector fields on pre-symplectic manifolds.

\section*{Conclusions}
\label{Sec:Conclusions}
\addcontentsline{toc}{section}{\nameref{Sec:Conclusions}}

In this paper we carried on the development of an approach to the inverse problem of the calculus of variations based on the use of the coisotropic embedding theorem that we initiated in \cite{Schiavone2023-Inverse_problem_Electrodynamics}.

\noindent Indeed, in \cite{Schiavone2023-Inverse_problem_Electrodynamics} we demonstrated how to provide a Lagrangian formulation of Maxwell equations in vacuum taking advantage of their formulation in terms of a vector field on a pre-symplectic manifold and by using the coisotropic embedding theorem.

\noindent In the present paper we provided the general theory underlying that example providing a solution for the inverse problem of the calculus of variations for all those (implicit) differential equations associated to a vector field on a pre-symplectic manifold with suitable conditions.

\bibliographystyle{alpha}
\bibliography{Biblio}

\end{document}